\documentclass[reqno]{amsart}
\usepackage{amssymb,mathtools,mathrsfs,anysize}
\marginsize{2.5cm}{2.5cm}{2.5cm}{2.5cm}

\usepackage{ifpdf}
\ifpdf
\usepackage[hyperindex]{hyperref}%,pagebackref
\else
\expandafter\ifx\csname dvipdfm\endcsname\relax
\usepackage[hypertex,hyperindex]{hyperref}
\else
\usepackage[dvipdfm,hyperindex]{hyperref}
\fi
\fi

\allowdisplaybreaks[4]

\theoremstyle{remark}
\newtheorem{rem}{Remark}

\DeclareMathOperator{\td}{d\!}
\DeclareMathOperator{\te}{e}

\begin{document}

\title[A unified proof of three combinatorial identities]
{A Unified Proof of Three Combinatorial Identities Related to the Stirling Numbers of the Second Kind}

\author[C.-Y. He]{Chun-Ying He}
\address{School of Mathematics and Physics, Hulunbuir University, Hailar 021008, Inner Mongolia, China}
\email{hechunying9209@qq.com}
\urladdr{\url{https://orcid.org/0000-0002-9709-8002}}

\author[F. Qi]{Feng Qi}
\address{Retired Professor and PhD, 17709 Sabal Court, University Village, Dallas, TX 75252-8024, USA}
\email{\href{mailto: F. Qi <qifeng618@gmail.com>}{qifeng618@gmail.com}}
\urladdr{\url{https://orcid.org/0000-0001-6239-2968}}

\dedicatory{Dedicated to Professor Huan-Nan Shi in honor of his 77th birthday}

\begin{abstract}
In the note, the authors give a unified proof of Identities~67, 84, and~85 in the monograph ``M. Z. Spivey, \emph{The Art of Proving Binomial Identities}, Discrete Mathematics and its Applications, CRC Press, Boca Raton, FL, 2019; available online at \url{https://doi.org/10.1201/9781351215824}'' and connect these three identities with a computing formula for the Stirling numbers of the second kind. Moreover, in terms of the notion of Qi's normalized remainders of the exponential and logarithmic functions, the authors reformulate the definitions of the Stirling numbers of the first and second kind and their generalizations by Howard in 1967 and 1980, Carlitz in 1980, and Broder in 1984.
\end{abstract}



\thanks{This paper has been accepted for publication in the Journal of Inequalities and Mathematical Analysis}

\thanks{This paper was typeset using \AmS-\LaTeX}

\maketitle

\section{Introduction}

Let $m,n\in\mathbb{N}_0$ be nonnegative integers.
\par
Identity~12 on~\cite[p.~12]{Spivey-art-2019} is
\begin{equation}\label{ID12-Spivey}
\sum _{k=0}^n\binom{n}{k} (-1)^k=
\begin{cases}
1, & n=0;\\
0, & n\in\mathbb{N}.
\end{cases}
\end{equation}
Identity~84 on~\cite[p.~62]{Spivey-art-2019} states
\begin{equation}\label{ID84-Spivey}
\sum_{k=0}^{n}\binom{n}{k}(n-k)^m(-1)^k=0, \quad m<n.
\end{equation}
Identity~67 on~\cite[p.~55]{Spivey-art-2019} reads that
\begin{equation}\label{ID67-Spivey}
\sum_{k=0}^{n}\binom{n}{k}(n-k)^n(-1)^k=n!.
\end{equation}
Identity~85 on~\cite[p.~62]{Spivey-art-2019} gives
\begin{equation}\label{ID85-Spivey}
\sum_{k=0}^{n}\binom{n}{k}(n-k)^{n+1}(-1)^k=\frac{n(n+1)!}{2}.
\end{equation}
It is easy to see that the second case $n\in\mathbb{N}$ in~\eqref{ID12-Spivey} is a special case of~\eqref{ID84-Spivey}. These identities can also be found in the monographs~\cite{Quaintance-Gould-Stirling-B, Sprugnoli-Gould-2006}.
\par
The only main aim of this note is to provide a unified proof for the above three identities~\eqref{ID84-Spivey} to~\eqref{ID85-Spivey}.

\section{Original proof of the identity~\eqref{ID67-Spivey}}\label{Foure-Otisg-p-sec}

We recite the original proof on~\cite[p.~68]{Spivey-art-2019} of the identity~\eqref{ID67-Spivey} as follows.
\par
By the binomial theorem, we have
\begin{equation}\label{Spivey-p68-first-eq}
\sum_{k=0}^{n}\binom{n}{k}\te^{(n-k)x}(-1)^k=(\te^x-1)^n.
\end{equation}
Now, differentiate both sides $n$ times. The left side is
\begin{equation*}
\sum_{k=0}^{n}\binom{n}{k}(n-k)^n\te^{(n-k)x}(-1)^k.
\end{equation*}
The right side is (after a derivative or two)
\begin{equation*}
\frac{\operatorname{d}^{n-1}}{\td x^{n-1}}\bigl[n(\te^x-1)^{n-1}\te^x\bigr]
=\frac{\operatorname{d}^{n-2}}{\td x^{n-2}}\bigl[n(n-1)(\te^x-1)^{n-2}\te^{2x} +n(\te^x-1)^{n-1}\te^x\bigr].
\end{equation*}
After $n$ derivatives we will have $n$ terms. One of these will be the term $n!\te^{nx}$. All of the others will contain at least one factor of $\te^x-1$. Now, let $x=0$. The $\te^{(n-k)x}$ factor vanishes on the left side, and everything on the right goes to sero except for $n!\te^0$. The identity follows.

\section{A unified proof of three identities}\label{Foure-unif-p-sec}

It is well known~\cite[p.~51]{Comtet-Combinatorics-74} that the Stirling numbers of the second kind $S(k,n)$ for $k\ge n\in\mathbb{N}_0$ can be analytically generated by
\begin{equation}\label{2stirl-gen-f}
\frac{(\te^x-1)^n}{n!}=\sum_{k=n}^{\infty} S(k,n)\frac{x^k}{k!}.
\end{equation}
Hence, we obtain
\begin{equation}\label{(tex-1)n-ser}
(\te^x-1)^n=\sum_{k=0}^{\infty} \frac{S(k+n,n)}{\binom{k+n}{n}}\frac{x^{k+n}}{k!}, \quad n\in\mathbb{N}_0.
\end{equation}
Substituting the Maclaurin power series expansion~\eqref{(tex-1)n-ser} into~\eqref{Spivey-p68-first-eq} leads to
\begin{equation}\label{(tex-1)n-ser-subst}
\sum_{k=0}^{n}\binom{n}{k}\te^{(n-k)x}(-1)^k=\sum_{k=0}^{\infty} \frac{S(k+n,n)}{\binom{k+n}{n}}\frac{x^{k+n}}{k!}, \quad n\in\mathbb{N}_0.
\end{equation}
Differentiating $m$ times with respect to $x$ on both sides of~\eqref{(tex-1)n-ser-subst} results in
\begin{equation}\label{(te)n-ser-subst}
\sum_{k=0}^{n}\binom{n}{k}(n-k)^m\te^{(n-k)x}(-1)^k=\sum_{k=0}^{\infty} \frac{S(k+n,n)}{\binom{k+n}{n}}\langle k+n\rangle_m \frac{x^{k+n-m}}{k!}
\end{equation}
for $m,n\in\mathbb{N}_0$, where
\begin{equation*}%\label{Fall-Factorial-Dfn-Eq}
\langle z\rangle_n=
\prod_{k=0}^{n-1}(z-k)=
\begin{cases}
z(z-1)\dotsm(z-n+1), & n\in\mathbb{N}\\
1,& n=0
\end{cases}
\end{equation*}
is the falling factorial of $z\in\mathbb{C}$.
\par
Taking the limit $x\to0$ on both sides of~\eqref{(te)n-ser-subst} reveals
\begin{equation}\label{te-lim-n-ser-subst}
\begin{split}
\sum_{k=0}^{n}\binom{n}{k}(n-k)^m(-1)^k
&=\begin{dcases}
\frac{S(m,n)}{\binom{m}{n}}\frac{\langle m\rangle_m}{(m-n)!}, & m\ge n\in\mathbb{N}_0\\
0, & n>m\in\mathbb{N}_0
\end{dcases}\\
&=\begin{dcases}
n!S(m,n), & m\ge n\in\mathbb{N}_0\\
0, & n>m\in\mathbb{N}_0
\end{dcases}
\end{split}
\end{equation}
for $m,n\in\mathbb{N}_0$. 
\par
The case $n>m\in\mathbb{N}_0$ in~\eqref{te-lim-n-ser-subst} is just the identity~\eqref{ID84-Spivey}.
\par
Letting $m=n$ in~\eqref{te-lim-n-ser-subst} gives the identity~\eqref{ID67-Spivey}.
\par
Taking $m=n+1$ in~\eqref{te-lim-n-ser-subst} and employing $S(n+1,n)=\frac{n(n+1)}{2}$ for $n\in\mathbb{N}_0$, we easily derive the identity~\eqref{ID85-Spivey}. The unified proof is complete.

\section{Remarks}

Finally, we list several remarks.

\begin{rem}
The unified proof in Section~\ref{Foure-unif-p-sec} is more understandable and comprehensive than the original proof recited in Section~\ref{Foure-Otisg-p-sec} of the identity~\eqref{ID67-Spivey}.
\end{rem}

\begin{rem}
We can write the equality~\eqref{te-lim-n-ser-subst} as the formula
\begin{equation}\label{te-lim-n-ser-write}
S(m,n)=
\begin{dcases}
\frac{1}{n!}\sum_{k=0}^{n}(-1)^k\binom{n}{k}(n-k)^m, & m\ge n\in\mathbb{N}_0;\\
0, & n>m\in\mathbb{N}_0,
\end{dcases}
\end{equation}
see~\cite[p.~193, Identity~224]{Spivey-art-2019}.
Consequently, since $S(n,n)=1$ for $n\in\mathbb{N}_0$, the identity~\eqref{ID67-Spivey} is the special case $m=n$ of the formula~\eqref{te-lim-n-ser-write}, which can be rearranged as the form
\begin{equation}\label{Stirling-Number-dfn}
S(m,n)=
\begin{dcases}
\frac{(-1)^{n}}{n!}\sum_{k=0}^n(-1)^{k}\binom{n}{k}k^{m}, & m>n\in\mathbb{N}_0;\\
1, & m=n\in\mathbb{N}_0;\\
0, & n>m\in\mathbb{N}_0,
\end{dcases}
\end{equation}
see~\cite[p.~204, Theorem~A]{Comtet-Combinatorics-74}.
\end{rem}

\begin{rem}
The three identities~\eqref{ID84-Spivey}, \eqref{ID67-Spivey}, and~\eqref{ID85-Spivey} can be reformulated as
\begin{align*}
\sum_{k=0}^{n}\binom{n}{k}k^m(-1)^k&=0, \quad m<n,\\
\sum_{k=0}^{n}\binom{n}{k}k^n(-1)^k&=(-1)^nn!,
\end{align*}
and
\begin{equation*}
\sum_{k=0}^{n}\binom{n}{k}k^{n+1}(-1)^k=(-1)^n\frac{n(n+1)!}{2}.
\end{equation*}
\end{rem}

\begin{rem}
The Maclaurin power series expansion~\eqref{2stirl-gen-f} can be reformulated as
\begin{equation}\label{2Stirl-funct-rew}
\biggl(\frac{\te^x-1}{x}\biggr)^n=\sum_{k=0}^\infty \frac{S(k+n,n)}{\binom{k+n}{n}} \frac{x^{k}}{k!}, \quad n\ge0.
\end{equation}
The Bernoulli numbers $B_{k}$ for $k\ge0$ are generated by
\begin{equation}\label{Bernumber-dfn}
\frac{x}{\te^x-1}=\sum_{k=0}^\infty B_k\frac{x^k}{k!}=1-\frac{x}2+\sum_{k=1}^\infty B_{2k}\frac{x^{2k}}{(2k)!}, \quad |x|<2\pi,
\end{equation}
see~\cite[p.~48]{Comtet-Combinatorics-74}. Comparing the generating functions in~\eqref{2Stirl-funct-rew} and~\eqref{Bernumber-dfn}, considering the fact that $\frac{\te^x-1}{x}$ and $\frac{x}{\te^x-1}$ are the reciprocal of each other, we are sure that the Bernoulli numbers $B_k$ and the Stirling numbers of the second kind $S(n,k)$ must have something to do with each other. This idea has been carried out and verified by Qi and his coauthors in the papers~\cite{dema-D-22-00032.tex, Genocchi-Stirling.tex, Guo-Qi-JANT-Bernoulli.tex, exp-derivative-sum-Combined.tex, Bernoulli-Stirling-4P.tex, Uniform-treatments}, for example.
\par
On the other hand, the generating functions $\bigl(\frac{\te^x-1}{x}\bigr)^n$ and $\frac{x}{\te^x-1}$ of the Stirling numbers of the second kind $S(n,k)$ and the Bernoulli numbers $B_k$ have been generalized in~\cite{Howard-Duke-1967-599, Howard-Fib-1980i4, Howard-Duke-1967-701} by
\begin{equation}\label{T-Associate-Stirling}
(T_r[\te^x])^\ell
=\frac{\ell![(r+1)!]^\ell}{[(r+1)\ell]!} \sum_{j=0}^{\infty}\frac{S_r(j+(r+1)\ell,\ell)}{\binom{j+(r+1)\ell}{j}} \frac{x^{j}}{j!}
\end{equation}
and
\begin{equation*}%\label{Howard-Gen-polyn-writ}
\frac{\te^{x t}}{T_{r-1}[\te^x]}
=\sum_{j=0}^{\infty}A_{r,j}(t)\frac{x^j}{j!}
\end{equation*}
respectively, where
\begin{equation}\label{expon-remainder-norm}
T_r[\te^x]=\frac{(r+1)!}{x^{r+1}}\Biggl(\te^x-\sum_{j=0}^{r}\frac{x^j}{j!}\Biggr), \quad r\in\mathbb{N}_0
\end{equation}
is called Qi's normalized remainder of the exponential function $\te^x$ in the literature~\cite{llog-cosine-expan.tex, log-secant-norm-tail.tex}.
For more information on Qi's normalized remainder $T_r[\te^z]$ of the exponential function $\te^x$, please refer to~\cite{log-exp-expan-Sym.tex, very-short-exp-nt.tex, exp-norm-tail-ratio.tex, exp-tail-ratio.tex}, \cite[Section~1]{ratio-NR-squ-tan.tex}, \cite[Section~1.7]{arcsine-liu-qi-mdpi.tex}, \cite[Remark~2]{JMI-5046.tex}, and closely related references therein.
\end{rem}

\begin{rem}
In~\cite{Howard-Fib-1980i4}, Howard defined $s_r(j,\ell)$ by
\begin{equation}\label{1st-Stirl-gen-log}
\Biggl[\ln\frac{1}{1-x}-\sum_{j=1}^{r}\frac{x^j}{j}\Biggr]^\ell
=\ell!\sum_{j=(r+1)\ell}^{\infty}s_r(j,\ell)\frac{x^j}{j!}, \quad r\in\mathbb{N}_0.
\end{equation}
It is clear that $s_0(j,\ell)=(-1)^{j+\ell}s(j,\ell)$, where the Stirling numbers of the first kind $s(j,\ell)$ can be analytically generated~\cite[Theorem~3.14]{Mansour-Schork-B2016} by
\begin{equation}\label{Stirl-No-First-GF}
\biggl[\frac{\ln(1+x)}{x}\biggr]^\ell=\sum_{j=0}^\infty \frac{s(j+\ell,\ell)}{\binom{j+\ell}{\ell}}\frac{x^{j}}{j!}, \quad |x|<1.
\end{equation}
\par
The equation~\eqref{1st-Stirl-gen-log} can be reformulated as
\begin{equation}\label{1st-Stirl-gen-log-RF}
\Biggl[(-1)^r\frac{r+1}{x^{r}}\Biggl(\frac{\ln(1+x)}{x}-\sum_{j=0}^{r-1}(-1)^j\frac{x^{j}}{j+1}\Biggr)\Biggr]^\ell
=\frac{\ell!(r+1)^\ell}{[(r+1)\ell]!} \sum_{j=0}^{\infty}(-1)^j\frac{s_r(j+(r+1)\ell,\ell)} {\binom{j+(r+1)\ell}{j}}\frac{x^{j}}{j!}
\end{equation}
for $r\in\mathbb{N}_0$. The function
\begin{equation}\label{logn-remainder-norm}
(-1)^{r}\frac{r+1}{x^{r}}\Biggl[\frac{\ln(1+x)}{x}-\sum_{j=0}^{r-1}(-1)^{j}\frac{x^{j}}{j+1}\Biggr]
\end{equation}
in~\eqref{1st-Stirl-gen-log-RF} is just Qi's normalized remainder $T_r\bigr[\frac{\ln(1+x)}{x}\bigr]$ for $r\in\mathbb{N}_0$ of the function $\frac{\ln(1+x)}{x}$.
\end{rem}

\begin{rem}
Theorem~15 in~\cite{Broder-1984} reads that the $r$-Stirling numbers of the first kind for $r\in\mathbb{N}_0$ have the ``vertical'' exponential generating function
\begin{equation}\label{r-stirling-first}
\sum_{k}\begin{bmatrix}
k+r\\m+r
\end{bmatrix}_r\frac{z^k}{k!}=
\begin{dcases}
\frac{1}{m!}\biggl(\frac{1}{1-z}\biggr)^r\biggl(\ln\frac{1}{1-z}\biggr)^m, &m\ge0;\\
0, & m<0.
\end{dcases}
\end{equation}
We can rewrite~\eqref{r-stirling-first} in the form
\begin{equation}\label{r-stirling-first-reform}
\biggl(\frac{1}{1+z}\biggr)^r\biggl[\frac{\ln(1+z)}{z}\biggr]^m
=\biggl(\frac{1}{1+z}\biggr)^r\biggl(T_0\biggr[\frac{\ln(1+z)}{z}\biggr]\biggr)^m
=\sum_{k=0}^\infty(-1)^k\frac{\left[\begin{smallmatrix}
k+m+r\\m+r
\end{smallmatrix}\right]_r}{\binom{k+m}{m}}\frac{z^k}{k!}, \quad |z|<1
\end{equation}
for $r,m\in\mathbb{N}_0$.
Taking $r=0$ in~\eqref{1st-Stirl-gen-log-RF} and~\eqref{r-stirling-first-reform} and comparing with~\eqref{Stirl-No-First-GF} give
\begin{equation*}
\begin{bmatrix}
k\\m
\end{bmatrix}_0=s_0(k,m)=(-1)^{k+m}s(k,m), \quad k,m\in\mathbb{N}_0.
\end{equation*}
\par
Basing on the above discussion, we propose a problem: Investigate the properties of the sequence $F(r,s,m,k)$ generated by
\begin{equation}\label{problem-log-normalized}
\biggl(\frac{1}{1+z}\biggr)^r\biggl(T_s\biggr[\frac{\ln(1+z)}{z}\biggr]\biggr)^m
=\sum_{k=0}^\infty F(r,s,m,k)\frac{z^k}{k!}, \quad r,m\in\mathbb{C}, \quad s\in\mathbb{N}_0, \quad |z|<1.
\end{equation}
\end{rem}

\begin{rem}
Theorem~16 in~\cite{Broder-1984} and Eq.~(3.9) in~\cite{Carlitz-Fibonacci-1980-I} state that the $r$-Stirling numbers of the second kind for $r\in\mathbb{N}_0$ have the exponential generating function
\begin{equation}\label{r-stirling-second}
\sum_{k}\left\{\begin{matrix}
k+r\\ m+r
\end{matrix}\right\}_r\frac{z^k}{k!}=
\begin{dcases}
\frac{1}{m!}\te^{rz}(\te^z-1)^m, & m\ge0;\\
0, & m<0.
\end{dcases}
\end{equation}
We can reformulate~\eqref{r-stirling-second} in the form
\begin{equation}\label{small-matrix-stirl}
\te^{rz}(T_0[\te^z])^m=\sum_{k=0}^\infty\frac{\left\{\begin{smallmatrix}
k+m+r\\ m+r
\end{smallmatrix}\right\}_r}{\binom{k+m}{m}}\frac{z^k}{k!}, \quad r,m\in\mathbb{N}_0, \quad |z|<\infty.
\end{equation}
Taking $r=0$ in~\eqref{T-Associate-Stirling} and~\eqref{small-matrix-stirl} and comparing with~\eqref{2Stirl-funct-rew} yield
\begin{equation*}
\left\{\begin{matrix}
k\\ m
\end{matrix}\right\}_0
=S_0(k,m)=S(k,m), \quad k,m\in\mathbb{N}_0.
\end{equation*}
\par
Basing on the above discussion, we propose a problem: Investigate the properties of the sequence $Q(r,s,m,k)$ generated by
\begin{equation*}
\te^{rz}(T_s[\te^z])^m=\sum_{k=0}^{\infty}Q(r,s,m,k)\frac{z^k}{k!}, \quad r,m\in\mathbb{C}, \quad s\in\mathbb{N}_0?
\end{equation*}
\end{rem}

\begin{rem}
We recall from the papers~\cite{log-exp-expan-Sym.tex, very-short-exp-nt.tex, ratio-NR-squ-tan.tex, llog-cosine-expan.tex, exp-norm-tail-ratio.tex, arcsine-liu-qi-mdpi.tex, JMI-5046.tex, log-secant-norm-tail.tex, exp-tail-ratio.tex} that Qi's normalized remainder can be generally defined as follows.
\par
Let $f$ be a real infinitely differentiable function on an interval $I$ such that $0\in I\subseteq\mathbb{R}$.
If $f^{(n+1)}(0)\ne0$ for some $n\in\mathbb{N}_0$, then the function
\begin{equation}\label{Qi-remainder-dfn}
T_n[f(x)]=
\begin{dcases}
\frac{1}{f^{(n+1)}(0)}\frac{(n+1)!}{x^{n+1}}\Biggl[f(x)-\sum_{k=0}^{n} f^{(k)}(0)\frac{x^k}{k!}\Biggr], & x\ne0\\
1, & x=0
\end{dcases}
\end{equation}
for $x\in I$ is said to be the $n$th normalized remainder or the $n$th normalized tail of the Maclaurin expansion of the function $f$.
\par
Applying $f(x)$ in~\eqref{Qi-remainder-dfn} to $\te^x$ leads to the normalized remainder $T_n[\te^x]$ for $n\in\mathbb{N}_0$ defined by~\eqref{expon-remainder-norm}, replacing $f(x)$ by $\frac{\ln(1+x)}{x}$ in~\eqref{Qi-remainder-dfn} reduces to the normalized remainder $T_n\bigr[\frac{\ln(1+x)}{x}\bigr]$ for $n\in\mathbb{N}_0$ defined by~\eqref{logn-remainder-norm}.
\par
It is easy to verify that $T_n\bigr[\frac{\ln(1+x)}{x}\bigr]=T_n[\ln(1+x)]$ for $n\in\mathbb{N}_0$. Therefore, the equations~\eqref{1st-Stirl-gen-log-RF} and~\eqref{problem-log-normalized} can be reformulated respectively by
\begin{equation*}
(T_r[\ln(1+x)])^\ell
=\frac{\ell!(r+1)^\ell}{[(r+1)\ell]!} \sum_{j=0}^{\infty}(-1)^j\frac{s_r(j+(r+1)\ell,\ell)} {\binom{j+(r+1)\ell}{j}}\frac{x^{j}}{j!},\quad r,\ell\in\mathbb{N}_0, \quad |x|<1
\end{equation*}
and
\begin{equation*}
\biggl(\frac{1}{1+z}\biggr)^r(T_s[\ln(1+z)])^m
=\sum_{k=0}^\infty F(r,s,m,k)\frac{z^k}{k!}, \quad r,m\in\mathbb{C}, \quad s\in\mathbb{N}_0, \quad |z|<1.
\end{equation*}
\end{rem}

\section*{Acknowledgements}
The authors were partially supported by the Youth Project of Hulunbuir City for Basic Research and Applied Basic Research (Grant No.~GH2024020) and by the Natural Science Foundation of Inner Mongolia Autonomous Region (Grant No.~2025QN01041).
\par
The authors are thankful to anonymous referees for their valuable comments and helpful suggestions to the original version of this paper.

\end{document}